\documentclass[english]{article}
\usepackage[T1]{fontenc}
\usepackage[latin9]{inputenc}
\usepackage{color}
\usepackage{float}
\usepackage{amsmath}
\usepackage{setspace}
\usepackage[authoryear]{natbib}
\makeatletter
\@ifundefined{date}{}{\date{}}
\usepackage[font=footnotesize,labelfont=bf]{caption}

\AtBeginDocument{

}

\makeatother

\usepackage{babel}
\begin{document}
\begin{onehalfspace}

\title{\noindent \textbf{\textcolor{black}{Double-calibration estimators
accounting for under-coverage and nonresponse in socio-economic surveys }}}
\end{onehalfspace}

\maketitle
\begin{onehalfspace}
\textcolor{black}{Maria Michela Dickson$^{1}$, Giuseppe Espa$^{1}$,
Lorenzo Fattorini$^{2}$}
\end{onehalfspace}

~

\begin{singlespace}
\textcolor{black}{\footnotesize{}$^{1}$ Department of Economics and
Management, University of Trento.}{\footnotesize\par}

\textcolor{black}{\footnotesize{}$^{2}$ Department of Economics and
Statistics, University of Siena.}{\footnotesize\par}
\end{singlespace}
\begin{onehalfspace}

\subsection*{\textcolor{black}{Summary}}
\end{onehalfspace}

\begin{onehalfspace}
\textcolor{black}{Under-coverage and nonresponse problems are jointly
present in most socio-economic surveys. The purpose of this paper
is to propose a completely design-based estimation strategy that accounts
for both problems without resorting to models but simply performing
a two-step calibration. The first calibration exploits a set of auxiliary
variables only available for the units in the sampled population to
account for nonresponse. The second calibration exploits a different
set of auxiliary variables available for the whole population, to
account for under-coverage. The two calibrations are then unified
in a double-calibration estimator. Mean and variance of the estimator
are derived up to the first order of approximation. Conditions ensuring
approximate unbiasedness are derived and discussed. The strategy is
empirically checked by a simulation study performed on a set of artificial
populations. A case study is lead on Danish data coming from the European
Union Statistics on Income and Living Conditions survey. The strategy
proposed is flexible and suitable in most situations in which both
under-coverage and nonresponse are present.}
\end{onehalfspace}

~
\begin{onehalfspace}

\subsection*{\textcolor{black}{Keywords}}
\end{onehalfspace}

\begin{onehalfspace}
\textcolor{black}{Auxiliary variables; Cut-off sampling; Design-based
estimation; First-order Taylor series approximation; Simulation study. }
\end{onehalfspace}

\begin{onehalfspace}

\section{\textcolor{black}{Introduction}}
\end{onehalfspace}

\begin{doublespace}
\textcolor{black}{\citet[p. 8]{sarndalj} establish four requirements
to select a probability sample, which set the perimeter for the definition
of a sampling design under the randomization principle. One of them
requires that the procedure to select the sample ensures invariably
positive probabilities to enter the sample for all units in the population. }

\textcolor{black}{This requirement may not be suitable in some situations
such as in establishment surveys, such as the Economic Census conducted
by the U.S. Census Bureau, in which the population of businesses is
characterized by a highly skewed distribution in the survey variables
(\citealp{glasser1962complete}). In this case, different approaches
are widespread used, essentially based on the partition of population
into strata determined by several business characteristics (e.g. size),
and some strata are completely censused, some are sampled, and some
are neglected, basing on units features or on the chance to contact
them (\citealp{sigman1995selecting}). As happen in establishment
surveys conducted by the U.S. Bureau of Economic Analysis, very small
establishments are excluded a priori from the population to be sampled,
due to the costs in build and update a sampling frame for them, against
of an expected slight gain in efficiency of the estimators (see e.g.
\citealp{hidiroglou1986construction,de1999item,rivest2002generalization}).
These instances are known in literature as cut-off sampling (\citealp{knaub2008cutoff,benedetti2010framework,haziza2010sampling}).
A similar position can be seen in social surveys on households, such
as e.g. the Household Finance and Consumption Survey managed by the
European Central Bank, characterized by the missed observation of
population units considered ineligibles for the survey, scilicet dwelling
that are vacant, not habitable, with non-eligible members, etc., with
consequences on the estimation of living conditions and poverty rate
(\citealp{nicoletti2011estimating}). Owing to the aforementioned
under-coverage of the whole population, the classical Horvitz-Thompson
(HT) estimator is biased in these situations. Bias is usually corrected
in literature by means of model-based techniques (see, among others,
\citealp{kott2006using,haziza2010sampling}). Recently, a design-based
solution to under-coverage problems has been proposed by \citet{fattorini2018use}
adopting a calibration technique in which the weights originally attached
to each sample observation are modified in such a way to be able to
estimate the population totals of a set of auxiliary variables without
error. The rationale behind calibration is well known: if the calibrated
weights guess the population totals of the auxiliary variables without
errors, they should be suitable also for estimating the total of the
survey variable, providing a relationship existing between the survey
variable and the auxiliaries. Obviously, calibration is likely to
perform well in terms of precision under strong linear relationship. }

\textcolor{black}{Socio-economic surveys are also interested by unit
nonresponse, which are the more frequent as more sensitive are the
survey variables (e.g. sexual behavior, drug consumption, etc.). Even
if undesirable, nonresponse is a natural contingency in surveys, so
that damages on estimation and inference caused by them need to be
addressed (\citealp{groves2008impact}). This argument is crucial
in survey sampling theory and it is extensively treated in literature
(e.g. \citealp{brick2009nonresponse}). Widely applied methods include
post-stratification \citep{holt1979post} or, more recently, once
again model-based techniques including imputation and nonresponse
propensity weighting (\citealp{sarndal2005estimation,haziza2010effect}).
By means of imputation, nonresponse values are replaced by substitutes
and estimation is performed on the completed sample. Imputed values
are customarily obtained by means of a prediction model presuming
a relationship between the survey variable and a set of covariates
known for all the population units or for all the sampled units. In
accordance with the presumed model, commonly used techniques of imputation
are, for example, regression imputation, nearest neighbor imputation,
hot deck imputation, and multiple imputation. Nonresponse propensity
weighting assumes that each unit of the sampled population has a strictly
positive probability to respond. A model is then used to estimate
the probabilities of respondent units from the sample by connecting
these probabilities to auxiliary information by means of logistic
regression models (\citealp{chang2008using}). In addition to this
source of uncertainty, the requirement of positive response probability
seems tighten in socio-economic surveys, because there will always
be units that do not respond in any situation (e.g. homeless and geographically
mobile individuals and families). Alternatively, \citet{fattorini2013design}
attempt a complete design-based solution in which population values
and nonresponse are viewed as fixed characteristics. To this purpose,
they once again use the calibration technique, termed in the literature
as nonresponse calibration weighting by \citet{haziza2010effect}.
In this case, weights originally attached to each respondent observation
are modified in such a way to be able at estimating the population
totals of a set of auxiliary variables without error.}

\textcolor{black}{In most cases under-coverage and nonresponse problems
are jointly present in socio-economic surveys. Therefore, a general
indication in the treatment of both problems concern the use of any
available auxiliary information, even if some of them are not available
for all units of the population. In this paper, we build on the availability
of a set of auxiliary variables for the whole population while another
set is available only for the sampled portion. In establishments surveys,
for example, many financial information may be available only for
businesses of adequate size, such as corporations, while they may
be not for those small businesses excluded from sampling, such as
micro-enterprises. Moreover, owing to the recent data collection developments,
the additional information may arise from big data, e.g. data coming
from internet and telephone use, social networks, online purchases,
etc. }

\textcolor{black}{The purpose of this paper is to propose double-calibration
estimators. The use of calibration in two or more steps has been used
by \citet{folsom2000generalized} and \citet{estevao2006survey}.
Here we propose a completely design-based estimation strategy that
considers both under-coverage and nonresponse problems, solving them
without resort to models but simply by performing a double calibration.
The first calibration exploits a set of auxiliary variables available
only for the units in the sampled population to account for nonresponse;
the second calibration exploits a different set of auxiliary variables
available for the whole population, to account for under-coverage.
Joining together the two calibrations, we propose a double-calibration
estimator that is applicable to all cases in which both under-coverage
and nonresponse problems are present.}

\textcolor{black}{The paper is structured as follow. In Section 2,
some preliminaries and notations are given. Section 3 is devoted to
the costruction of the double-calibration estimator and in Section
4 some design-based properties (expectation and variance) are derived.
In order to check the efficiency of the strategy, in Section 5 a Monte
Carlo simulation study exploring several scenarios is performed. In
Section 6, by using data coming from the European Union Statistics
on Income and Living Conditions survey and from Statistics Denmark
data, a case study to estimate the total income of Danish households
in 2013 is presented and discussed. Some concluding remarks are given
in Section 7.}
\end{doublespace}

\textcolor{black}{~}

\section{Preliminaries and notation}

\begin{doublespace}
\textcolor{black}{Denote as $U=\left\{ u_{1},...,u_{N}\right\} $
a finite population of $N$ units. Let $y_{j}$, with $j\in U$, the
value for unit $j$ of the survey variable $Y$. We aim at estimating
the population total $T_{Y}=\sum_{j\in U}y_{j}$. For the whole population
there exists a vector $\boldsymbol{Z}$ of $M$ auxiliary variables
whose values $\boldsymbol{\boldsymbol{z}}_{j}=\left[z_{j1},...,z_{jM}\right]^{t}$
are known for each $j\in U$, in such a way that the vector of totals
$\boldsymbol{T}_{Z}=\sum_{j\in U}\boldsymbol{z}_{j}$ is also known. }

\textcolor{black}{In this setting, for one of the reasons mentioned
in the introduction, only a sub-population $U_{B}$ of size $N_{B}<N$
units is sampled using a fixed-size design having first- and second-order
inclusion probabilities $\pi_{j},\pi_{jh}$ for any $h>j\in U_{B}$.
Denote by $T_{Y(B)}=\sum_{j\in U_{B}}y_{j}$ the unknown total of
$Y$ in $U_{B}$. Moreover, suppose that additional information exists
in the sub-population $U_{B}$, possibly arising from big data sources.
More precisely suppose that there exists a vector $\boldsymbol{X}$
of $K$ auxiliary variables whose values $\boldsymbol{x}_{j}=\left[x_{j1},...,x_{jM}\right]^{t}$
are known for each $j\in U_{B}$ in such a way that the vector of
totals $\boldsymbol{T}_{X(B)}=\sum_{j\in U_{B}}\boldsymbol{x}_{j}$
is also known. In this setting, denote by $\boldsymbol{T}_{Z(B)}=\sum_{j\in U_{B}}\boldsymbol{\boldsymbol{z}}_{j}$
the known vector of total of the $\boldsymbol{z}_{j}$s in the sub-population
$U_{B}$. }

\textcolor{black}{A random sample $S$ of $n<N_{B}$ units is selected
from the sub-population $U_{B}$ by means of the adopted sampling
scheme. As often happens in practice, especially in socio-economic
surveys, the sample may be affected by nonresponses, in such a way
that the sample is split into two sub-samples, the sub-sample $R\subset S$
of the respondent units and the sub-sample $S-R$ of the nonrespondent
units. }

\textcolor{black}{The above presented setup shows two problems to
solve: first, a correction for nonresponses is necessary, in order
to estimate $T_{Y(B)}$; second, since the sample $S$ is selected
from $U_{B}$ and not from $U$, any $T_{Y(B)}$ estimator is biased,
so that it needs for a correction in order to estimate $T_{Y}$. We
propose a design-based calibration in two steps, developed in next
sub-sections. }
\end{doublespace}
\begin{onehalfspace}

\section{The double-calibration estimator}
\end{onehalfspace}
\begin{onehalfspace}

\subsection{\textit{\textcolor{black}{First calibration: from respondent group
to sampled sub-population}}}
\end{onehalfspace}

\begin{doublespace}
\textcolor{black}{The first issue to deal with is the nonresponse
problem occurring in a sample. Since $S$ is selected in $U_{B}$,
in absence of nonresponses, it would be possible to estimate $T_{Y(B)}$
by means of the well-known HT estimator}

\begin{equation}
\hat{T}_{Y(B)}=\sum_{j\in S}\frac{y_{j}}{\pi_{j}}
\end{equation}

\textcolor{black}{and $\hat{T}_{Y(B)}$ would be an unbiased estimator
for $T_{Y(B)}$. However, owing to nonresponses, $\hat{T}_{Y(B)}$
is unknown and the response-based estimator }

\[
\hat{T}_{Y(B)R}=\sum_{j\in R}\frac{y_{j}}{\pi_{j}}\neq\hat{T}_{Y(B)}
\]

\textcolor{black}{is a biased estimator of $T_{Y(B)}$. Following
results obtained in \citet{sarndal2005estimation}, the bias may be
reduced by exploiting the}\textbf{\textcolor{black}{$\boldsymbol{X}$}}\textcolor{black}{-vector
of auxiliary information. The resulting estimator is}

\begin{equation}
\hat{T}_{Y(B)cal}=\hat{\boldsymbol{b}}_{R}^{t}\boldsymbol{T}_{X(B)}
\end{equation}

\textcolor{black}{where $\hat{\boldsymbol{b}}_{R}=\hat{\boldsymbol{A}}_{R}^{-1}\hat{\boldsymbol{a}}_{R}$
is the least-square coefficient vector of the regression of $Y$ vs}\textbf{\textcolor{black}{{}
$\boldsymbol{X},$}}\textcolor{black}{{} performed on the respondent
sample $R$, i.e. $\hat{\boldsymbol{A}}_{R}=\sum_{j\in R}\frac{\boldsymbol{x}_{j}\boldsymbol{x}_{j}^{t}}{\pi_{j}}$
and $\hat{\boldsymbol{a}}_{R}=\sum_{j\in R}\frac{y_{j}\boldsymbol{x}_{j}}{\pi_{j}}$
and the unit constant is tacitly adopted as the first auxiliary variable
in the vector}\textbf{\textcolor{black}{$\boldsymbol{X}$}}\textcolor{black}{. }

\textcolor{black}{The design-based properties of $\hat{T}_{Y(B)cal}$
are derived in \citet{fattorini2013design}. The population is partitioned
into respondent and nonrespondent strata and the estimator is approximately
unbiased if the relationship between $Y$ and }\textbf{\textcolor{black}{$\boldsymbol{X}$}}\textcolor{black}{{}
is similar in both the strata. Practically speaking, this condition
is similar to those assumed in most model-based nonresponse treatment
even if not embedded into models. }
\end{doublespace}

~
\begin{onehalfspace}

\subsection{\textit{\textcolor{black}{Second calibration: from sampled sub-population
to the whole population}}}
\end{onehalfspace}

\begin{doublespace}
\textcolor{black}{Because $\hat{T}_{Y(B)cal}$ is, at most, an approximately
unbiased estimator of $T_{Y(B)}$, it is a biased estimator of $T_{Y}$.
Indeed, the sampling scheme adopted to select $S$ generates a sampling
design onto $U_{B}$ but not onto $U$, and units of $U-U_{B}$ cannot
enter the sample. Therefore, the missed selection of some population
units leads to a bias due to population under-coverage and it is necessary
to correct the estimator $\hat{T}_{Y(B)cal}$. }

\textcolor{black}{\citet{fattorini2018use} named these schemes as
pseudo designs and propose a design-based calibration estimation based
on a single auxiliary variable having a proportional relationship
with the survey variable. In order to extend this approach to vectors
of auxiliary variables and to more general linear relationships, the
population under-coverage is handled by the calibration criterion
proposed by \citet{sarndal2005estimation}. More precisely, if the
$y_{j}$s were available for each $j\in S$, the information furnished
by the $M$ auxiliary variables $\boldsymbol{Z}$, available for the
whole population units, could be exploited by means of the calibration
estimator }

\begin{equation}
\hat{T}_{Y(cal)}=\hat{\boldsymbol{d}}_{B}^{t}\boldsymbol{T}_{Z}
\end{equation}

\textcolor{black}{where $\hat{\boldsymbol{d}_{B}}=\hat{\boldsymbol{C}}_{B}^{-1}\hat{\boldsymbol{c}}_{B}$
is the least-square coefficient vector of the regression of $Y$ vs
$\boldsymbol{Z}$, performed on the whole sample $S$, i.e. $\boldsymbol{\hat{C}}_{B}=\sum_{j\in S}\frac{\boldsymbol{\boldsymbol{z}}_{j}\boldsymbol{\boldsymbol{z}}_{j}^{t}}{\pi_{j}}$
and $\boldsymbol{\hat{c}}_{B}=\sum_{j\in S}\frac{y_{j}\boldsymbol{\boldsymbol{z}}_{j}}{\pi_{j}}$.}

\textcolor{black}{If we suppose once again that the unit constant
is adopted as the first auxiliary variable in the vector $\boldsymbol{Z}$,
then the calibration estimator (3) could be rewritten as}

\begin{equation}
\hat{T}_{Y(cal)}=\hat{T}_{Y(B)}+\hat{\boldsymbol{d}}_{B}^{t}(\boldsymbol{T}_{Z}-\hat{\boldsymbol{T}}_{Z(B)})
\end{equation}

\textcolor{black}{wher}e $\hat{\boldsymbol{T}}_{Z(B)}=\sum_{j\in S}\frac{\boldsymbol{\boldsymbol{z}}_{j}}{\pi_{j}}$\textcolor{black}{{}
is the HT estimator of the totals of the $\boldsymbol{z}_{j}$s in
the sampled sub-population $U_{B}$ (see Appendix A.1 for the proof). }

\textcolor{black}{However, the estimator }$\hat{T}_{Y(cal)}$\textcolor{black}{{}
is only virtual, because having the values of the survey variable
only for the respondent subset $R$, neither the HT estimator }$\hat{T}_{Y(B)}$\textcolor{black}{{}
nor the least-squares coefficient vector $\hat{\boldsymbol{d}_{B}}=\hat{\boldsymbol{C}}_{B}^{-1}\hat{\boldsymbol{c}}_{B}$
are known. Therefore, exploiting equation (4), a }\textit{\textcolor{black}{double
calibration estimator}}\textcolor{black}{{} can be constructed by using
}$\hat{T}_{Y(B)cal}$\textcolor{black}{{} instead of }$\hat{T}_{Y(B)}$\textcolor{black}{{}
and $\hat{\boldsymbol{d}_{R}}=\hat{\boldsymbol{C}}_{R}^{-1}\hat{\boldsymbol{c}}_{R}$,
instead of $\hat{\boldsymbol{d}_{B}}$ where $\boldsymbol{\hat{C}}_{R}=\sum_{j\in R}\frac{\boldsymbol{\boldsymbol{z}}_{j}\boldsymbol{\boldsymbol{z}}_{j}^{t}}{\pi_{j}}$
and $\boldsymbol{\hat{c}}_{R}=\sum_{j\in R}\frac{y_{j}\boldsymbol{\boldsymbol{z}}_{j}}{\pi_{j}}$.
Practically speaking, the resulting estimator of the whole population
total turns out to be }

\begin{equation}
\hat{T}_{Y(dcal)}=\hat{T}_{Y(B)cal}+\hat{\boldsymbol{d}}_{R}^{t}\mathrm{(}\boldsymbol{T}_{Z}-\hat{\boldsymbol{T}}_{Z(B)})=\hat{\boldsymbol{b}}_{R}^{t}\boldsymbol{T}_{X(B)}+\hat{\boldsymbol{d}}_{R}^{t}\mathrm{(}\boldsymbol{T}_{Z}-\hat{\boldsymbol{T}}_{Z(B)})
\end{equation}

\textcolor{black}{By the double calibration estimator, the information
provided by $\boldsymbol{X}$ and $\boldsymbol{Z}$ is exploited for
handling both nonresponses and population under-coverage in a design-based
framework. }
\end{doublespace}

\noindent ~
\begin{onehalfspace}

\section{Design-based properties of the double calibration estimator}
\end{onehalfspace}

\begin{doublespace}
\textcolor{black}{Let denote by $U_{B(R)}$ the stratum of respondent
units in the sub-population $U_{B}$ and by $U_{B(NR)}$ the stratum
of nonrespondent units. As suggested by \citet{fattorini2013design},
let introduce a dummy variable as $r_{j}=1$ if $j\in U_{B(R)}$ and
$r_{j}=0$ if $j\in U_{B(NR)}$. Therefore, using the $r_{j}$s indicators
$\hat{\boldsymbol{A}}_{R}$, $\hat{\boldsymbol{a}}_{R}$, $\hat{\boldsymbol{C}}_{R}$
and $\hat{\boldsymbol{c}}_{R}$ can be rewritten as $\hat{\boldsymbol{A}}_{R}=\sum_{j\in S}\frac{r_{j}\boldsymbol{x}_{j}\boldsymbol{x}_{j}^{t}}{\pi_{j}}$,
$\hat{\boldsymbol{a}}_{R}=\sum_{j\in S}\frac{r_{j}y_{j}\boldsymbol{\boldsymbol{x}}_{j}}{\pi_{j}}$,$\boldsymbol{\hat{C}}_{R}=\sum_{j\in S}\frac{r_{j}\boldsymbol{\boldsymbol{z}}_{j}\boldsymbol{\boldsymbol{z}}_{j}^{t}}{\pi_{j}}$
and $\boldsymbol{\hat{c}}_{R}=\sum_{j\in S}\frac{r_{j}y_{j}\boldsymbol{\boldsymbol{z}}_{j}}{\pi_{j}}$.
Due to that, the previous matrices and vectors as well as the double
calibration estimator }$\hat{T}_{Y(dcal)}$\textcolor{black}{{} depend
on the selection of the sole sample $S$, while nonresponses are accounted
for in the $r_{j}$s, which are a fixed characteristic of the population,
as is required in a design-based perspective. }

\textcolor{black}{It is worth noting that in this perspective, $\hat{\boldsymbol{A}}_{R}$,
$\hat{\boldsymbol{a}}_{R}$, $\hat{\boldsymbol{C}}_{R}$, $\hat{\boldsymbol{c}}_{R}$
and }$\hat{\boldsymbol{T}}_{Z(B)}$\textcolor{black}{{} are HT estimators
of $\boldsymbol{A}_{R}=\sum_{j\in U_{B}}r_{j}\boldsymbol{x}_{j}\boldsymbol{x}_{j}^{t}=\sum_{j\in U_{B(R)}}\boldsymbol{x}_{j}\boldsymbol{x}_{j}^{t}$,
$\boldsymbol{a}_{R}=\sum_{j\in U_{B}}r_{j}y_{j}\boldsymbol{x}_{j}=\sum_{j\in U_{B(R)}}y_{j}\boldsymbol{x}_{j}$,
$\boldsymbol{C}_{R}=\sum_{j\in U_{B}}r_{j}\boldsymbol{z}_{j}\boldsymbol{z}_{j}^{t}=\sum_{j\in U_{B(R)}}\boldsymbol{z}_{j}\boldsymbol{z}_{j}^{t}$,
$\boldsymbol{c}_{R}=\sum_{j\in U_{B}}r_{j}y_{j}\boldsymbol{z}_{j}=\sum_{j\in U_{B(R)}}y_{j}\boldsymbol{z}_{j}$
and of }$\boldsymbol{T}_{Z(B)}$\textcolor{black}{, respectively.
Therefore, because} $\hat{T}_{Y(dcal)}$\textcolor{black}{{} is differentiable
with respect to $\hat{\boldsymbol{A}}_{R}$, $\boldsymbol{\hat{a}}_{R}$,
$\boldsymbol{\hat{C}}_{R}$, $\hat{\boldsymbol{c}}_{R}$ and }$\hat{\boldsymbol{T}}_{Z(B)}$\textcolor{black}{,
it can be approximated up to the first term by a Taylor series around
the true population counterparts $\boldsymbol{A}_{R}$, $\boldsymbol{a}_{R}$,
$\boldsymbol{C}_{R}$, $\boldsymbol{c}_{R}$ and }$\boldsymbol{T}_{Z(B)}$\textcolor{black}{.
The equation of the first-order Taylor series approximation of} $\hat{T}_{Y(dcal)}$\textcolor{black}{{}
is derived in Appendix A.2. }
\end{doublespace}

\noindent ~
\begin{onehalfspace}

\subsection{Approximate expectation}
\end{onehalfspace}

\begin{doublespace}
From the first-order Taylor series approximation of $\hat{T}_{Y(dcal)}$
it immediately follows that 

\begin{equation}
AE(\hat{T}_{Y(dcal)})=\boldsymbol{b}_{R}^{t}\boldsymbol{T}_{X(B)}+\boldsymbol{d}_{R}^{t}(\boldsymbol{T}_{Z}-\boldsymbol{T}_{Z(B)})
\end{equation}

where $\boldsymbol{b}_{R}=\boldsymbol{A}_{R}^{-1}\boldsymbol{a}_{R}$
is the least-square coefficient vector of the regression of $Y$ vs
$\boldsymbol{X}$ performed on the respondent stratum $U_{B(R)}$
and $\boldsymbol{d}_{R}=\boldsymbol{C}_{R}^{-1}\boldsymbol{c}_{R}$
is the least-square coefficient vector of the regression of $Y$ vs
$\boldsymbol{Z}$ performed in the same stratum. Exploiting equation
(6), after some algebra reported in Appendix A.3, it is proven that
the double calibration estimator is unbiased up to the first-order
approximation if: 
\end{doublespace}
\begin{enumerate}
\begin{doublespace}
\item the linear relationship between $Y$ and $\boldsymbol{X}$ is similar
in the respondent and nonrespondent strata of $U_{B}$, i.e. $\boldsymbol{b}_{R}\approx\boldsymbol{b}_{NR}$,
where $\boldsymbol{b}_{NR}$ is the least-square coefficient vector
of the regression of $Y$ vs $\boldsymbol{X}$ performed on the nonrespondent
stratum $U_{B(NR)}$ ; 
\item the linear relationship between $Y$ and $\boldsymbol{Z}$ is similar
in the respondent stratum and in the whole sub-population $U_{B}$,
i.e. $\boldsymbol{d}_{R}\approx\boldsymbol{d}_{B}$, where $\boldsymbol{d}_{B}$
is the least-square coefficient vector of the regression of $Y$ vs
$\boldsymbol{Z}$ performed on the whole sub-population $U_{B}$; 
\item the linear relationship between $Y$ and $\boldsymbol{Z}$ is similar
in the two sub-populations $U_{B}$ and $U-U_{B}$, i.e. $\boldsymbol{d}_{B}\approx\boldsymbol{d}_{NB}$,
where $\boldsymbol{d}_{NB}$ is the least-square coefficient vector
of the regression of $Y$ vs $\boldsymbol{Z}$ performed on the whole
sub-population $U-U_{B}$. 
\end{doublespace}
\end{enumerate}
\noindent ~
\begin{onehalfspace}

\subsection{Approximate variance and variance estimation}
\end{onehalfspace}

\begin{doublespace}
From equation (A.3) of Appendix A.2, the first-order Taylor series
approximation of $\hat{T}_{Y(dcal)}$ is rewritten as a translation
of a HT estimator, in the sense that

\[
\hat{T}_{Y(dcal)}=cost+\sum_{j\in S}\frac{u_{j}}{\pi_{j}}
\]

where 

\begin{multline*}
u_{j}=r_{j}\left(y_{j}\boldsymbol{x}_{j}^{t}-\boldsymbol{a}_{R}^{t}\boldsymbol{A}_{R}^{-1}\boldsymbol{x}_{j}\boldsymbol{x}_{j}^{t}\right)\boldsymbol{A}_{R}^{-1}\boldsymbol{T}_{X(B)}+\\
+r_{j}\left(y_{j}\boldsymbol{z}_{j}^{t}-\boldsymbol{c}_{R}^{t}\boldsymbol{C}_{R}^{-1}\boldsymbol{z}_{j}\boldsymbol{z}_{j}^{t}\right)\boldsymbol{C}_{R}^{-1}\left(\boldsymbol{T}_{Z}-\boldsymbol{T}_{Z(B)}\right)-\boldsymbol{c}_{R}^{t}\boldsymbol{C}_{R}^{-1}\boldsymbol{z}_{j},j\in U_{B}
\end{multline*}

are the influence values (e.g. \citealp{davison1997bootstrap}).

Therefore, the approximate variance of $\hat{T}_{Y(dcal)}$ turns
out to be (e.g. \citealp[p. 175]{sarndalj}) 

\begin{equation}
AV\left(\hat{T}_{Y(dcal)}\right)=\sum_{h>j\in U_{B}}(\pi_{j}\pi_{h}-\pi_{jh})\left(\frac{u_{j}}{\pi_{j}}-\frac{u_{h}}{\pi_{h}}\right)^{2}
\end{equation}

On the basis of equation (7), the well-known Sen-Yates-Grundy (SYG)
variance estimator is given by 

\begin{equation}
\hat{V}_{SYG}^{2}=\sum_{h>j\in S}(\pi_{j}\pi_{h}-\pi_{jh})\left(\frac{\hat{u}_{j}}{\pi_{j}}-\frac{\hat{u}_{h}}{\pi_{h}}\right)^{2}
\end{equation}

where

\begin{multline*}
\hat{u}_{j}=r_{j}\left(y_{j}\boldsymbol{x}_{j}^{t}-\hat{\boldsymbol{a}}_{R}^{t}\boldsymbol{\hat{A}}_{R}^{-1}\boldsymbol{x}_{j}\boldsymbol{x}_{j}^{t}\right)\boldsymbol{\hat{A}}_{R}^{-1}\boldsymbol{T}_{X(B)}+\\
+r_{j}\left(y_{j}\boldsymbol{z}_{j}^{t}-\boldsymbol{\hat{c}}_{R}^{t}\boldsymbol{\hat{C}}_{R}^{-1}\boldsymbol{z}_{j}\boldsymbol{z}_{j}^{t}\right)\boldsymbol{\hat{C}}_{R}^{-1}\left(\boldsymbol{T}_{Z}-\hat{\boldsymbol{T}}_{Z(B)}\right)-\boldsymbol{\hat{c}}_{R}^{t}\hat{\boldsymbol{C}}_{R}^{-1}\boldsymbol{z}_{j},j\in S
\end{multline*}

are the empirical influence values computed for each sample unit. 
\end{doublespace}

\textcolor{black}{~}

\section{Simulation study}

\begin{doublespace}
Simulations were used to check the performance of the proposed estimator.
We considered a population $U$ of $N=10,000$ units and a sub-population
$U_{B}\subset U$ of $N_{B}=7,500$ units. We supposed that the values
$z_{j}$ of an auxiliary variable $Z$ were available for each $j\in U$
and were adopted for sample under-coverage calibration. Moreover,
we supposed that the values $x_{j}$ of an auxiliary variable $X$
achieved from additional information were available for each $j\in U_{B}$
and were adopted in nonresponse calibration. We also supposed that
the sub-population $U_{B}$ was partitioned into respondent and non-respondent
strata $U_{B(R)}$ and $U_{B(NR)}$ , respectively. Three sizes were
supposed for the respondent stratum, $N_{B(R)}=2,250;4,500;6,750$
units corresponding to response rates of 30\%, 60\% and 90\%, respectively.
The auxiliary variables $X$ and $Z$ and the survey variables $Y$
were generated from a tri-variate normal distribution. The expectations
and variances of $X$ and $Z$ were supposed to be equal to 1, while
the expectation and variance of $Y$ were supposed to be equal to
2 and 4, respectively. These set ups assured that each variable had
a coefficient of variation of 1. The correlation between $X$ and
$Y$ was set to be $\rho_{XY}=0.3;0.6;0.9$; similarly, the correlation
between $Z$ and $Y$ was set to be $\rho_{ZY}=0.3;0.6;0.9$, giving
rise to nine scenarios. The correlation between $X$ and $Z$ was
set to the minimum possible value $\rho_{XZ}$ such that the resulting
variance-covariance matrix is positive\textcolor{black}{{} definite.
Once the nine variance-covariance matrices were established the 10,000
values of $Z$ and $Y$ and the 7,500 values of $X$ were generated
using the triangular square root of the variance-covariance matrix
(e.g. \citealp[Sect. 4.1]{johnson2013multivariate}).} Subsequently,
the first $N_{B(R)}$ units of $U_{B}$ were supposed to be the respondent
portion of the population, ensuring in this way conditions 1.-3.,
i.e. the approximate unbiasedness of the double calibration estimator.
Simple random sampling without replacement (SRSWOR) was the sampling
scheme adopted to select samples of sizes $n=75;100;150;250;375;500$
from $U_{B}$. If the same sampling efforts were adopted to select
samples from the whole population $U$ and in absence of nonresponses,
then the HT estimator of the total would give rise to relative root
means squared errors

\begin{equation}
RRMSE_{SRSWOR}=\sqrt{\frac{N-n}{Nn}}CV_{Y}
\end{equation}

where $CV_{Y}$ is the coefficient of variation of the survey variable.
Equation (9) was taken as benchmark for the performance of the double
calibration estimator. 

For each combination of respondent sizes $N_{B(R)}$, correlations
between $X$ and $Y$, correlations between $Z$ and $Y$, and sample
sizes $n$, 10,000 random samples were selected by means of SRSWOR
from $U_{B}$, and the double calibration estimates $\hat{T}_{i}=\left(i=1,...,10000\right)$
were computed using equation (5). Moreover, from each simulated sample,
the variance estimates $\hat{V}_{i}^{2}=\left(i=1,...,10000\right)$
were also computed using equation (8), that under SRSWOR reduces to

\begin{equation}
\hat{V}_{SYG}^{2}=N_{B}\left(N_{B}-n\right)\frac{s_{\hat{u}}^{2}}{n}
\end{equation}

where $s_{\hat{u}}^{2}$ is the sampling variance of the $\hat{u}_{j}$s.
Once the variance estimates were computed from (10), the $RRMSE$
estimates $\hat{RRMSE}_{i}=\frac{\hat{V}_{i}}{\hat{T}_{i}}$ were
achieved together with the confidence intervals at the nominal level
of 0.95, $\hat{T}_{i}\pm2\hat{V}_{i}$. Therefore, from the resulting
Monte Carlo distributions of these quantities, the expectations $E(\hat{T}_{Y(dcal)})=\frac{1}{10000}\sum_{i=1}^{10000}\hat{T}_{i}$
and mean squared errors $MSE(\hat{T}_{Y(dcal)})=\frac{1}{10000}\sum_{i=1}^{10000}(\hat{T}_{i}-T_{Y})^{2}$
of the double calibration estimator were empirically derived from
which the relative bias $RB=\frac{E(\hat{T}_{Y(dcal)})-T_{Y}}{T_{Y}}$
and the relative root mean squared errors $RRMSE=\frac{\sqrt{MSE(\hat{T}_{Y(dcal)})}}{T_{Y}}$
were derived. The expectations of the $RRMSE$ estimator $ERRMSEE=\frac{1}{10000}\sum_{i=1}^{10000}\hat{RRMSE}_{i}$
and the coverage of the 0.95 confidence interval $COV95=\frac{1}{10000}\sum_{i=1}^{10000}I(\hat{T}_{i}-2\hat{V}_{i}\leq T_{Y}\leq\hat{T}_{i}+2\hat{V}_{i})$
are also computed. Simulation results are reported in Tables B.1-B.9
of the Appendix B. 

The simulation results motivate the following remarks. 

The first order approximation of relative bias and $RRMSE$ are very
accurate in most cases. The discrepancies between approximation and
the empirical values achieved from the Monte Carlo distributions are
usually smaller than one percent point that become lower with high
levels of response and correlations. In these cases, approximations
turn out to be smaller than the true values at most by three percent
points (Table B.9). The theoretical findings about the bias reduction,
reported in section 4.1 are fully confirmed by the simulation results.
The artificial populations considered in the study meet the unbiasedness
conditions 1.-3. Indeed the empirical values of the relative bias
are negligible (invariably about one percentage p\textcolor{black}{oint)
irres}pective of the level of correlation of the survey variable with
the auxiliaries. While the level of correlation does not affect the
bias reduction, it has a relevant impact on the precision. When correlations
are strong the double calibration estimator proves to be efficient,
reaching values of $RRMSE$ that are even smaller to those achieved
by the HT estimator adopting the same sampling effort and in absence
of nonresponse and under-coverage. Obviously precision increases also
with the level of response. 

The $RRMSE$ estimator obtained from the variance estimator (8) is
approximately unbiased providing also confidence intervals with coverage
near to the nominal level of 95\% in most cases. Because the estimator
(8) actually estimates the approximate variance, some exceptions occur
when the variance approximations (and the $RRMSE$ subsequently) turn
out be smaller than the true values. In these cases, an under-evaluation
of about three percent points produced a coverage only about 80\%
(see Table B.9). 
\end{doublespace}

~

\section{An application to the European Union Statistics on Income and Living
Conditions survey}

\begin{doublespace}
National statistical institutes periodically collect data on living
conditions through household surveys. Information contents concern
several aspects about living conditions, such as, among others, features
and expenses incurred to manage the dwelling, material deprivation
and welfare indicators, individual and household incomes. The European
Union Statistics on Income and Living Conditions survey was built
on previous experience of the European Community Household Panel (ECHP).
The survey was launched in 2003 in seven countries (Belgium, Denmark,
Greece, Ireland, Luxembourg, Austria and Norway), reaching all 28-EU
member countries, plus Switzerland, Norway, Iceland, FYROM and Serbia.
It is conducted yearly and collect information about European households.
Some rules to conduct the survey are set by the Eurostat, as, among
others, the frequency and the period to which questions must be referred,
and the aggregation level of some longitudinal and cross-sectional
estimates. Other aspects of the survey are set independently by each
country, such as, for instance, the sampling design and the sample
size, leading to several discrepancies among countries (see, among
others, \citealp{goedeme2013much}; \citealp{lohmann2011comparability}). 

Moreover, the population coverage of surveys like these is incomplete.
Individuals who do not live in households, as well as homeless, physically
or mentally unable, geographically mobile and displaced individuals
are not always represented in national-level data. It is estimated
that worldwide some 300 to 350 million people may be missing from
survey sampling frames, at least 45\% omitted altogether by design,
or because they are likely to be undercounted (\citealp{carr2013missing}).
The European Union Statistics on Income and Living Conditions survey,
which involves approximately 300,000 households across Europe, is
not an exception and is affected by under-coverage, as well as samples
selected are affected by nonresponses. We propose an example about
the use of the double-calibration estimator in the 2013 wave of the
European Union Statistics on Income and Living Conditions survey in
Denmark (hereinafter DK-SILC). Data about respondents are freely available
from the Eurostat web site, while other needed information were tracked
among those collected by Statistics Denmark. 

The reference population $U$ consists of households residing in Denmark,
except for those habitually living in a foreign country or institutional
cohabitations as orphanages, religious institutes, etc. As available
in Statistics Denmark website, the households population size in 2013
was equal to 2,891,119 units. The DK-SILC survey is based on a simple
random sampling without replacement design, so that inclusion probabilities
are equal for all units in the population. The sampling unit is the
individual person and the household is defined as the household of
which the selected person is member. This because an household in
Denmark is defined as composed of one or more individuals. Households
eligible to DK-SILC are those in which the sampling unit is a person
aged 16 and over, living alone or together in private dwellings and
bound by marriage, parentage, affinity or other relationships. So
that, the eligible population $U_{B}$ of Danish households equal
to 2,416,597, leading to an under-coverage rate equal to 0.16\%. 

The 2013 DK-SILC survey was also interested by a nonresponse rate
of about 63\%. In fact, the respondent number was equal to 5,419,
against a sample of 14,702 households. Micro-data about respondents
include many information, grouped in four sections: Household Register
(D), Personal Register (R), Household Data (H) and Personal Data (P).
Variables collected concern several information, most of which are
qualitative. To implement the present case study, we use quantitative
variables (in euro) referred to the previous year of survey (2012),
contained in the H-section. Specifically, the tax on income and social
contributions (HY140G) are used as the $X$ variable to correct for
nonresponse, while the total housing cost (HH070) are used as the
$Z$ variable to correct for under-coverage. The variable $Y$ under
estimation is the total household disposable income (HY020). Sample
data suggest that both auxiliary variables are slightly correlated
with the variable under estimation (0.38 among $X$ and $Y$; 0.17
among $Z$ and $Y$, in the respondent group), revealing an unfavorable
situation, worse than all presented in section 5. However, from simulation
results, the weak relationships between the survey and the auxiliary
variable should deteriorate precision but, fortunately, should not
deteriorate bias reduction. The estimated total household disposable
income is equal to 125,739.17 million euro, equivalent to an average
household disposable income on $U$ equal to 43,491.52 euro. Since
the sampling design is SRSWOR, the variance estimate is computed as
in (10) and the $RRMSE$ estimate results equal to 0.05. 

Results obtained need to be intended as an illustration and do not
pretend to be official estimates. Clearly, quality of the results
relies on the quality of available data. Howsoever, results are in
line with those disseminated by Statistics Denmark. In fact, the average
disposable income for all households (population $U$) in 2012 is
329,803 Danish krone, corresponding to approximately 44,203.67 euro
(at current exchange rate). 
\end{doublespace}

~
\begin{onehalfspace}

\section{\textcolor{black}{Final remarks}}
\end{onehalfspace}

\begin{doublespace}
\textcolor{black}{The double calibration estimator can be adopted
in socio-economic surveys to jointly account for nonresponse and under-coverage
in a complete design-based framework, without resorting to model but
simply adopting a two-step calibration. The first calibration, performed
to reduce nonresponse bias, needs for a set of auxiliary variables
whose totals are known for the sampled sub-populations and whose values
are known for the respondent units in the sample. The second calibration,
performed to reduce the bias generated by the cut-off sampling, needs
for a further set of auxiliary variables whose totals are known for
the whole populations and whose values are known for all the units
in the sample. Interestingly, no list frame is necessary for the non-sampled
sub-population. If the relationships of the survey variable with the
two sets of auxiliaries are approximately similar in sampled and non-sampled
sub-populations as well as in respondent and nonrespondent strata
(conditions 1.-3.), the proposed estimator proves to be effective
for reducing bias, being also efficient for high-quality auxiliary
variables correlated with the variable of interest. Socio-economic
surveys may benefit from application of double-calibration estimator.
It leads to results very near with those disseminated by national
institutes of statistics and typically collected by integrating several
data sources, with far less effort in terms of data collection and
integration.}
\end{doublespace}

\textcolor{black}{\pagebreak}

\textcolor{black}{}
\textcolor{black}{\bibliographystyle{apalike}
\bibliography{prova1}
}
\end{document}


\begin{singlespace}
\noindent \textbf{\textcolor{black}{\Large{}Supplementary material
for ``Double-calibration estimators accounting for under-coverage
and nonresponse in socio-economic surveys'' }}{\Large\par}
\end{singlespace}
\begin{center}
\textcolor{black}{Maria Michela Dickson$^{1}$, Giuseppe Espa$^{1}$,
Lorenzo Fattorini$^{2}$}
\par\end{center}

\textcolor{black}{~}

\begin{singlespace}
\textcolor{black}{\footnotesize{}$^{1}$ Department of Economics and
Management, University of Trento.}{\footnotesize\par}

\textcolor{black}{\footnotesize{}$^{2}$ Department of Economics and
Statistics, University of Siena.}{\footnotesize\par}
\end{singlespace}

\appendix

\section{Appendix A - Results of Sections 3 and 4}

\subsection{An alternative formulation of the calibration estimator (3)}

\begin{doublespace}
Because the unit constant is adopted as the first auxiliary variable,
then $\boldsymbol{e}_{1}^{t}\boldsymbol{z}_{j}=1,\forall j\in U$,
where $\boldsymbol{e}_{1}^{t}=\left[1,0,...,0\right]$ is the first
vector of the standard basis of\textcolor{black}{{} $\mathbb{R}^{M}$.
T}herefore, the difference between equations (3) and (4) is given
by 

\begin{multline*}
\hat{T}_{Y(B)}-\hat{\boldsymbol{d}}_{B}^{t}\boldsymbol{\hat{T}}_{Z(B)}=\sum_{j\in S}\frac{y_{j}}{\pi_{j}}-\boldsymbol{\hat{d}}_{B}^{t}\sum_{j\in S}\frac{\boldsymbol{z}_{j}}{\pi_{j}}=\\
=\sum_{j\in S}\frac{y_{j}}{\pi_{j}}-\left(\sum_{j\in S}\frac{\boldsymbol{z}_{j}^{t}}{\pi_{j}}\right)\hat{\boldsymbol{d}}_{B}=\sum_{j\in S}\frac{y_{j}}{\pi_{j}}-\left(\sum_{j\in S}\frac{\boldsymbol{e}_{1}^{t}\boldsymbol{z}_{j}\boldsymbol{z}_{j}^{t}}{\pi_{j}}\right)\hat{\boldsymbol{d}}_{B}=\\
=\sum_{j\in S}\frac{y_{j}}{\pi_{j}}-\boldsymbol{e}_{1}^{t}\left(\sum_{j\in S}\frac{\boldsymbol{z}_{j}\boldsymbol{z}_{j}^{t}}{\pi_{j}}\right)\hat{\boldsymbol{C}}_{B}^{-1}\hat{\boldsymbol{c}}_{B}=\sum_{j\in S}\frac{y_{j}}{\pi_{j}}-\boldsymbol{e}_{1}^{t}\hat{\boldsymbol{C}}_{B}\hat{\boldsymbol{C}}_{B}^{-1}\hat{\boldsymbol{c}}_{B}=\\
=\sum_{j\in S}\frac{y_{j}}{\pi_{j}}-\boldsymbol{e}_{1}^{t}\hat{\boldsymbol{c}}_{B}=\sum_{j\in S}\frac{y_{j}}{\pi_{j}}-\boldsymbol{e}_{1}^{t}\sum_{j\in S}\frac{y_{j}\boldsymbol{z}_{j}}{\pi_{j}}=\sum_{j\in S}\frac{y_{j}}{\pi_{j}}-\sum_{j\in S}\frac{y_{j}\boldsymbol{e}_{1}^{t}\boldsymbol{z}_{j}}{\pi_{j}}=\\
=\sum_{j\in S}\frac{y_{j}}{\pi_{j}}-\sum_{j\in S}\frac{y_{j}}{\pi_{j}}=0
\end{multline*}

Hence, (3) and (4) coincide. 

~
\end{doublespace}

\subsection{First-order Taylor series approximation of $\hat{T}_{Y(dcal)}$}

\begin{doublespace}
The double calibration estimator (5) can be rewritten as 

\begin{multline}
\hat{T}_{Y(dcal)}=f\left(\hat{\boldsymbol{a}}_{R},\hat{\boldsymbol{A}}_{R},\hat{\boldsymbol{c}}_{R},\hat{\boldsymbol{C}}_{R},\hat{\boldsymbol{T}}_{Z(B)}\right)=\hat{\boldsymbol{a}}_{R}^{t}\hat{\boldsymbol{A}}_{R}^{-1}\boldsymbol{T}_{X(B)}+\hat{\boldsymbol{c}}_{R}^{t}\hat{\boldsymbol{C}}_{R}^{-1}\left(\boldsymbol{T}_{Z}-\hat{\boldsymbol{T}}_{Z(B)}\right)=\\
\sum_{k=1}^{K}\sum_{k\prime=1}^{K}\hat{a}_{Rk}T_{X(B)k\prime}\hat{a}_{R}^{kk\prime}+\sum_{m=1}^{M}\sum_{m\prime=1}^{M}\hat{c}_{Rm}(T_{Zm\prime}-\hat{T}_{Z(B)m\prime}){\hat{c}}_{R}^{mm\prime}\label{eq:A.1}
\end{multline}

\textcolor{black}{where ${\hat{a}}_{R}^{kk\prime}$ denotes the $kk$'
element of $\hat{\boldsymbol{A}}_{R}^{-1}$, $\hat{a}_{Rk}$ and $T_{X(B)k}$
are the $k$ components of $\hat{\boldsymbol{a}}_{R}$ and $\boldsymbol{T}_{X(B)}$,
respectively, ${\color{red}{\color{black}c}_{{\color{black}R}}^{{\color{black}mm\prime}}}$denotes
the $mm$' element of $\hat{\boldsymbol{C}}_{R}^{-1}$, and $\hat{c}_{Rm}$,
$T_{Zm}$} \textcolor{black}{and $\hat{T}_{Z(B)m}$ are }the $m$
components of $\hat{\boldsymbol{c}}_{R}$, $\boldsymbol{T}_{Z}$ and
$\hat{\boldsymbol{T}}_{Z(B)}$, respectively. Differentiating \ref{eq:A.1}
with respect to all the variables involved, it follows that 

~

$\begin{alignedat}{2}\frac{\partial f}{\partial\hat{a}_{Rk}}=\sum_{k\prime=1}^{K}T_{X(B)k\prime}\hat{a}_{R}^{kk\prime}, &  &  & k=1,...,K\end{alignedat}
$

$\begin{alignedat}{2}\frac{\partial f}{\partial{\hat{a}}_{Rkk\prime}}=-\hat{\boldsymbol{a}}_{R}^{t}\hat{\boldsymbol{A}}_{R}^{-1}\boldsymbol{E}_{kk\prime}\hat{\boldsymbol{A}}_{R}^{-1}\boldsymbol{T}_{X(B)}, &  &  & k,k^{\prime}=1,\ldots,M\end{alignedat}
$

$\begin{alignedat}{2}\frac{\partial f}{\partial{\hat{c}}_{Rm}}=\sum_{m^{\prime}=1}^{K}{\left(T_{Zm^{\prime}}-{\hat{T}}_{Z\left(B\right)m^{\prime}}\right){\hat{c}}_{R}^{mm^{\prime}}}, &  &  & m=1,\ldots,M\end{alignedat}
$

$\begin{alignedat}{2}\frac{\partial f}{\partial{\hat{a}}_{Rmm^{\prime}}}=-\hat{\boldsymbol{c}}_{R}^{t}\hat{\boldsymbol{C}}_{R}^{-1}\boldsymbol{E}_{mm^{\prime}}\hat{\boldsymbol{C}}_{R}^{-1}\left(\boldsymbol{T}_{Z}-{\hat{\boldsymbol{T}}}_{Z\left(B\right)}\right), &  &  & m,m^{\prime}=1,\ldots,M\end{alignedat}
$

$\begin{alignedat}{2}\frac{\partial f}{\partial{\hat{T}}_{Z(B)m\prime}}=-\sum_{m=1}^{K}{{\hat{c}}_{Rm}{\hat{c}}_{R}^{mm^{\prime}}}, &  &  & m\prime=1,\ldots,M\end{alignedat}
$

~

where $\boldsymbol{E}_{kk\prime}$ is the $K$-square matrix of $0$s,
with a 1 in position $kk$\textquoteright , and $\boldsymbol{E}_{mm\prime}$
is the $M$-square matrix of $0$s, with a 1 in position $mm\lyxmathsym{\textquoteright}$.
Evaluating these partial derivatives at the expected points ${\hat{\boldsymbol{a}}}_{R}=\boldsymbol{a}_{R},{\hat{\boldsymbol{A}}}_{R}=\boldsymbol{A}_{R},{\hat{\boldsymbol{c}}}_{R}=\boldsymbol{c}_{R},{\hat{\boldsymbol{C}}}_{R}=\boldsymbol{C}_{R}$
and ${\hat{\boldsymbol{T}}_{Z(B)}}=\boldsymbol{T}_{Z(B)}$, the first-order
Taylor series approximation of ${\hat{T}}_{Y(dcal)}$ gives rise to 

\begin{multline}
{\hat{T}}_{Y(dcal)}={\boldsymbol{a}}_{R}^{t}\boldsymbol{A}_{R}^{-1}\boldsymbol{T}_{X(B)}+\boldsymbol{c}_{R}^{t}\boldsymbol{C}_{R}^{-1}(\boldsymbol{T}_{Z}-\boldsymbol{T}_{Z\left(B\right)})+\sum_{k=1}^{K}{\sum_{k^{\prime}=1}^{K}{{(\hat{a}}_{Rk}-a_{Rk})}T_{X\left(B\right)k^{\prime}}a_{R}^{kk^{\prime}}}-\\
\sum_{k=1}^{K}{\sum_{k^{\prime}=1}^{K}{\boldsymbol{a}_{R}^{t}\boldsymbol{A}_{R}^{-1}\boldsymbol{E}_{kk^{\prime}}\boldsymbol{A}_{R}^{-1}\boldsymbol{T}_{X\left(B\right)}}\left({\hat{a}}_{R}^{kk^{\prime}}-a_{R}^{kk^{\prime}}\right)}+\\
\sum_{m=1}^{M}{\sum_{m^{\prime}=1}^{M}{{(\hat{c}}_{Rm}-c_{Rm})}(T_{Zm^{\prime}}-{T_{Z\left(B\right)m^{\prime}})c}_{R}^{mm^{\prime}}}-\\
\sum_{m=1}^{M}{\sum_{m^{\prime}=1}^{M}{\boldsymbol{c}_{R}^{t}\boldsymbol{C}_{R}^{-1}\boldsymbol{E}_{mm^{\prime}}\boldsymbol{C}_{R}^{-1}\left({\boldsymbol{T}_{Z}-\boldsymbol{T}}_{X\left(B\right)}\right)}\left({\hat{c}}_{R}^{mm^{\prime}}-c_{R}^{mm^{\prime}}\right)}-\\
\sum_{m=1}^{M}\sum_{m^{\prime}=1}^{M}{c_{Rm}c_{R}^{mm^{\prime}}}({\hat{T}}_{Z\left(B\right)k}-T_{Z\left(B\right)k})=\\
{\boldsymbol{a}}_{R}^{t}\boldsymbol{A}_{R}^{-1}\boldsymbol{T}_{X\left(B\right)}+\boldsymbol{c}_{R}^{t}\boldsymbol{C}_{R}^{-1}\left(\boldsymbol{T}_{Z}-\boldsymbol{T}_{Z\left(B\right)}\right){+\left({\hat{\boldsymbol{a}}}_{R}-\boldsymbol{a}_{R}\right)}^{t}\boldsymbol{A}_{R}^{-1}\boldsymbol{T}_{X\left(B\right)}-\\
\boldsymbol{a}_{R}^{t}\boldsymbol{A}_{R}^{-1}({\hat{\boldsymbol{A}}}_{R}-\boldsymbol{A}_{R})\boldsymbol{A}_{R}^{-1}\boldsymbol{T}_{X\left(B\right)}+{\left({\hat{\boldsymbol{c}}}_{R}-\boldsymbol{c}_{R}\right)}^{t}\boldsymbol{C}_{R}^{-1}{(\boldsymbol{T}}_{Z}-\boldsymbol{T}_{Z\left(B\right)})-\\
\boldsymbol{c}_{R}^{t}\boldsymbol{C}_{R}^{-1}({\hat{\boldsymbol{C}}}_{R}-\boldsymbol{C}_{R})\boldsymbol{C}_{R}^{-1}{(\boldsymbol{T}}_{Z}-\boldsymbol{T}_{Z\left(B\right)})-\boldsymbol{c}_{R}^{t}\boldsymbol{C}_{R}^{-1}{(\hat{\boldsymbol{T}}}_{Z(B)}-\boldsymbol{T}_{Z(B)})=\\
{\boldsymbol{a}}_{R}^{t}\boldsymbol{A}_{R}^{-1}\boldsymbol{T}_{X(B)}{+\hat{\boldsymbol{a}}}_{R}^{t}\boldsymbol{A}_{R}^{-1}\boldsymbol{T}_{X(B)}-\boldsymbol{a}_{R}^{t}\boldsymbol{A}_{R}^{-1}{\hat{\boldsymbol{A}}}_{R}\boldsymbol{A}_{R}^{-1}\boldsymbol{T}_{X(B)}+\\
{\hat{\boldsymbol{c}}}_{R}^{t}\boldsymbol{C}_{R}^{-1}{(\boldsymbol{T}}_{Z}-\boldsymbol{T}_{Z(B)})-\boldsymbol{c}_{R}^{t}\boldsymbol{C}_{R}^{-1}{\hat{\boldsymbol{C}}}_{R}\boldsymbol{C}_{R}^{-1}{(\boldsymbol{T}}_{Z}-\boldsymbol{T}_{Z(B)})+\\
\boldsymbol{c}_{R}^{t}\boldsymbol{C}_{R}^{-1}\boldsymbol{T}_{Z}-\boldsymbol{c}_{R}^{t}\boldsymbol{C}_{R}^{-1}{\hat{\boldsymbol{T}}}_{Z(B)}=c_{1}{+\hat{\boldsymbol{a}}}_{R}^{t}\boldsymbol{A}_{R}^{-1}\boldsymbol{T}_{X(B)}-\boldsymbol{a}_{R}^{t}\boldsymbol{A}_{R}^{-1}{\hat{\boldsymbol{A}}}_{R}\boldsymbol{A}_{R}^{-1}\boldsymbol{T}_{X(B)}+\\
c_{2}+{\hat{\boldsymbol{c}}}_{R}^{t}\boldsymbol{C}_{R}^{-1}{(\boldsymbol{T}}_{Z}-\boldsymbol{T}_{Z(B)})-\boldsymbol{c}_{R}^{t}\boldsymbol{C}_{R}^{-1}{\hat{\boldsymbol{C}}}_{R}\boldsymbol{C}_{R}^{-1}{(\boldsymbol{T}}_{Z}-\boldsymbol{T}_{Z(B)})-\boldsymbol{c}_{R}^{t}\boldsymbol{C}_{R}^{-1}{\hat{\boldsymbol{T}}}_{Z(B)}\label{eq:A.2}
\end{multline}

Grouping the constant terms, equation \ref{eq:A.2} can be rewritten
as 

\begin{multline}
{\hat{T}}_{Y(dcal)}=cost+\left(\sum_{j\in S}\frac{{r_{j}y}_{j}\boldsymbol{x}_{j}}{\pi_{j}}\right)^{t}\boldsymbol{A}_{R}^{-1}\boldsymbol{T}_{X(B)}-\\
\boldsymbol{a}_{R}^{t}\boldsymbol{A}_{R}^{-1}\left(\sum_{j\in S}\frac{{r_{j}\boldsymbol{x}}_{j}\boldsymbol{x}_{j}^{t}}{\pi_{j}}\right)\boldsymbol{A}_{R}^{-1}\boldsymbol{T}_{X(B)}+\left(\sum_{j\in S}\frac{{r_{j}y}_{j}\boldsymbol{z}_{j}}{\pi_{j}}\right)^{t}\boldsymbol{C}_{R}^{-1}{(\boldsymbol{T}}_{Z}-\boldsymbol{T}_{Z(B)})-\\
\boldsymbol{c}_{R}^{t}\boldsymbol{C}_{R}^{-1}\left(\sum_{j\in S}\frac{{r_{j}\boldsymbol{z}}_{j}\boldsymbol{z}_{j}^{t}}{\pi_{j}}\right)\boldsymbol{C}_{R}^{-1}{(\boldsymbol{T}}_{Z}-\boldsymbol{T}_{Z(B)})-\boldsymbol{c}_{R}^{t}\boldsymbol{C}_{R}^{-1}\left(\sum_{j\in S}\frac{\boldsymbol{z}_{j}}{\pi_{j}}\right)=\\
=cost+\sum_{j\in S}\frac{u_{j}}{\pi_{j}}
\end{multline}

where

\begin{eqnarray*}
u_{j}=r_{j}(y_{j}\boldsymbol{x}_{j}^{t}{-\boldsymbol{a}}_{R}^{t}\boldsymbol{A}_{R}^{-1}\boldsymbol{x}_{j}\boldsymbol{x}_{j}^{t}{)\boldsymbol{A}}_{R}^{-1}\boldsymbol{T}_{X\left(B\right)}+\\
r_{j}(y_{j}\boldsymbol{z}_{j}^{t}{-\boldsymbol{c}}_{R}^{t}\boldsymbol{C}_{R}^{-1}\boldsymbol{z}_{j}\boldsymbol{z}_{j}^{t}{)\boldsymbol{C}}_{R}^{-1}\left(\boldsymbol{T}_{Z}-\boldsymbol{T}_{Z\left(B\right)}\right)-\boldsymbol{c}_{R}^{t}\boldsymbol{C}_{R}^{-1}\boldsymbol{z}_{j}, &  & j\in U_{B}
\end{eqnarray*}

are the influence values (e.g. \citealp{davison1997bootstrap}).
\end{doublespace}

~

\subsection{Approximate unbiasedness of the double calibration estimator}

\begin{doublespace}
Regarding the first term of equation (6), \citet{fattorini2013design}
show that

\[
\boldsymbol{b}_{R}^{t}\boldsymbol{T}_{X(B)}=T_{Y(B)}+\left(\boldsymbol{b}_{N}-\boldsymbol{b}_{NR}\right)^{t}\boldsymbol{T}_{X(NR)}
\]

where $\boldsymbol{T}_{X(NR)}=\sum_{j\in U_{B(NR)}}\boldsymbol{x}_{j}$.
Therefore, if condition 1. holds, i.e. if $\boldsymbol{b}_{R}\approx\boldsymbol{b}_{NR}$,
then equation (6) approximately reduces to 

\[
AE(\hat{T}_{Y(dcal)})\approx T_{Y(B)}+\boldsymbol{d}_{R}^{t}(\boldsymbol{T}_{Z}-\boldsymbol{T}_{Z(B)})
\]

in such a way that, if condition 2. holds, i.e. if $\boldsymbol{d}_{R}\approx\boldsymbol{d}_{B}$
then 

\begin{equation}
AE({\hat{T}}_{Y\left(dcal\right)})\approx T_{Y(B)}\ +\boldsymbol{d}_{B}^{t}(\boldsymbol{T}_{Z}-\boldsymbol{T}_{Z(B)})\label{eq:A.4}
\end{equation}

Therefore, from \ref{eq:A.4} it follows that 

\begin{multline*}
AE\left({\hat{T}}_{Y\left(dcal\right)}\right)\approx\sum_{j\in U_{B}}y_{j}\ +\boldsymbol{d}_{B}^{t}\sum_{j\in U-U_{B}}\boldsymbol{z}_{j}=\sum_{j\in U_{B}}y_{j}\ +\\
\boldsymbol{d}_{B}^{t}\sum_{j\in U-U_{B}}\boldsymbol{z}_{j}+\sum_{j\in U-U_{B}}y_{j}-\sum_{j\in U-U_{B}}y_{j}=\sum_{j\in U}y_{j}\ +\boldsymbol{d}_{B}^{t}\sum_{j\in U-U_{B}}\boldsymbol{z}_{j}-\\
\sum_{j\in U-U_{B}}y_{j}=T_{Y}\ +\boldsymbol{d}_{B}^{t}\sum_{j\in U-U_{B}}\boldsymbol{z}_{j}-\sum_{j\in U-U_{B}}{(\boldsymbol{d}}_{NB}^{t}\boldsymbol{z}_{j}+e_{NBj})
\end{multline*}

where $\boldsymbol{d}_{NB}$ is the least-square coefficient vector
of the regression of $Y$ vs $\boldsymbol{Z}$ performed on the unsampled
portion of the population $U-U_{B}$ and the $e_{NB}j$s are the regression
residuals. 
\end{doublespace}

Because regression residuals sum to 0, then 

\begin{doublespace}
\begin{multline}
AE\left({\hat{T}}_{Y\left(dcal\right)}\right)\approx T_{Y}+\boldsymbol{d}_{B}^{t}\sum_{j\in U-U_{B}}\boldsymbol{z}_{j}-\boldsymbol{d}_{NB}^{t}\sum_{j\in U-U_{B}}\boldsymbol{z}_{j}=\\
{T_{Y}+\left(\boldsymbol{d}_{B}-\boldsymbol{d}_{NB}\right)}^{t}(\boldsymbol{T}_{Z}-\boldsymbol{T}_{Z\left(B\right)})\label{eq:A.5}
\end{multline}

From equation \ref{eq:A.5} it is proved that${\hat{T}}_{Y\left(dcal\right)}$
is approximately unbiased if condition 3. holds, i.e if $\boldsymbol{d}_{B}\approx\boldsymbol{d}_{NB}$

~
\end{doublespace}

\section{\textcolor{black}{\LARGE{}Appendix B - Results of simulation study
presented in Section 5 }}

\begin{table}[H]
\caption{Percentage values of relative bias ($RB$), first order approximation
of relative root mean squared error ($ARRMSE$), relative root mean
squared error ($RRMSE$), expectation of the relative root mean squared
error estimator ($ERRMSEE$) and coverage of the nominal 95\% confidence
intervals achieved from a population of 10,000 units, a sampled sub-population
of 7500 units with 2250; 4500 and 7500 respondent units, sample sizes
$n=75;100;150;250;375;500$ selected by means of simple random sampling
without replacement and correlations $\rho_{XY}=0.3$, $\rho_{ZY}=0.3$
and $\rho_{XZ}=0$. Values in parentheses are the $RRMSE$s of the
Horvitz-Thompson estimator in absence of nonresponse and under-coverage.}

~

\begin{tabular}{|cccccc|}
$N_{B(R)}=2,250$ &  &  &  &  & \multicolumn{1}{c}{}\tabularnewline
\cline{1-1} 
Approximate relative bias -1.7 &  &  &  &  & \multicolumn{1}{c}{}\tabularnewline
\cline{1-1} 
$n$ & $RB$ & $ARRSME$ & $RRMSE$ & $ERRMSEE$ & $COV95$\tabularnewline
\hline 
75 & -1.7 & 20.0 & 21.0 (11.5) & 20.7 & 92.1\tabularnewline
100 & -1.9 & 17.3 & 17.7 (9.9) & 17.8 & 93.5\tabularnewline
150 & -1.6 & 14.1 & 14.5 (8.1) & 14.4 & 93.6\tabularnewline
250 & -1.8 & 10.8 & 11.1 (6.2) & 11.1 & 94.4\tabularnewline
375 & -1.7 & 8.8 & 8.9 (5.1) & 8.9 & 94.5\tabularnewline
500 & -1.9 & 7.5 & 7.8 (4.4) & 7.7 & 94.4\tabularnewline
\end{tabular}

~

~

\begin{tabular}{|cccccc|}
$N_{B(R)}=4,500$ &  &  &  &  & \multicolumn{1}{c}{}\tabularnewline
\cline{1-1} 
Approximate relative bias -1.0 &  &  &  &  & \multicolumn{1}{c}{}\tabularnewline
\cline{1-1} 
$n$ & $RB$ & $ARRSME$ & $RRMSE$ & $ERRMSEE$ & $COV95$\tabularnewline
\hline 
75 & -1.0 & 14.1 & 14.4 (11.5) & 14.4 & 94.2\tabularnewline
100 & -0.9 & 12.2 & 12.2 (9.9) & 12.4 & 94.9\tabularnewline
150 & -1.1 & 9.9 & 10.1 (8.1) & 10.1 & 94.6\tabularnewline
250 & -1.1 & 7.6 & 7.7 (6.2) & 7.7 & 94.8\tabularnewline
375 & -1.0 & 6.2 & 6.3 (5.1) & 6.3 & 95.0\tabularnewline
500 & -1.1 & 5.3 & 5.4 (4.4) & 5.4 & 95.0\tabularnewline
\end{tabular}

~

~

\begin{tabular}{|cccccc|}
$N_{B(R)}=6,750$ &  &  &  &  & \multicolumn{1}{c}{}\tabularnewline
\cline{1-1} 
Approximate relative bias -1.5 &  &  &  &  & \multicolumn{1}{c}{}\tabularnewline
\cline{1-1} 
$n$ & $RB$ & $ARRSME$ & $RRMSE$ & $ERRMSEE$ & $COV95$\tabularnewline
\hline 
75 & -1.6 & 11.2 & 11.5 (11.5) & 11.5 & 94.8\tabularnewline
100 & -1.6 & 9.7 & 9.9 (9.9) & 9.9 & 94.8\tabularnewline
150 & -1.4 & 7.9 & 8.1 (8.1) & 8.0 & 94.8\tabularnewline
250 & -1.6 & 6.1 & 6.3 (6.2) & 6.2 & 94.4\tabularnewline
375 & -1.5 & 4.9 & 5.2 (5.1) & 5.0 & 94.2\tabularnewline
500 & -1.5 & 4.2 & 4.5 (4.4) & 4.3 & 94.1\tabularnewline
\end{tabular}
\end{table}

~

\begin{table}[H]
\caption{Percentage values of relative bias ($RB$), first order approximation
of relative root mean squared error ($ARRMSE$), relative root mean
squared error ($RRMSE$), expectation of the relative root mean squared
error estimator ($ERRMSEE$) and coverage of the nominal 95\% confidence
intervals achieved from a population of 10,000 units, a sampled sub-population
of 7,500 units with 2,250; 4,500 and 7,500 respondent units, sample
sizes $n=75;100;150;250;375;500$ selected by means of simple random
sampling without replacement and correlations $\rho_{XY}=0.3$, $\rho_{ZY}=0.6$
and $\rho_{XZ}=0$. Values in parentheses are the $RRMSE$s of the
Horvitz-Thompson estimator in absence of nonresponse and under-coverage.}

~

\begin{tabular}{|cccccc|}
$N_{B(R)}=2,250$ &  &  &  &  & \multicolumn{1}{c}{}\tabularnewline
\cline{1-1} 
Approximate relative bias -0.9 &  &  &  &  & \multicolumn{1}{c}{}\tabularnewline
\cline{1-1} 
$n$ & $RB$ & $ARRSME$ & $RRMSE$ & $ERRMSEE$ & $COV95$\tabularnewline
\hline 
75 & -0.8 & 17.9 & 18.8 (11.5) & 18.4 & 92.9\tabularnewline
100 & -1.0 & 15.5 & 15.8 (9.9) & 15.8 & 94.0\tabularnewline
150 & -0.8 & 12.6 & 13.0 (8.1) & 12.8 & 93.8\tabularnewline
250 & -0.8 & 9.7 & 9.8 (6.2) & 9.8 & 95.1\tabularnewline
375 & -0.8 & 7.9 & 7.9 (5.1) & 7.9 & 95.0\tabularnewline
500 & -1.0 & 6.7 & 6.9 (4.4) & 6.8 & 95.1\tabularnewline
\end{tabular}

~

~

\begin{tabular}{|cccccc|}
$N_{B(R)}=4,500$ &  &  &  &  & \multicolumn{1}{c}{}\tabularnewline
\cline{1-1} 
Approximate relative bias -1.0 &  &  &  &  & \multicolumn{1}{c}{}\tabularnewline
\cline{1-1} 
$n$ & $RB$ & $ARRSME$ & $RRMSE$ & $ERRMSEE$ & $COV95$\tabularnewline
\hline 
75 & -1.0 & 12.2 & 12.5 (11.5) & 12.5 & 94.1\tabularnewline
100 & -1.0 & 10.6 & 10.6 (9.9) & 10.7 & 94.8\tabularnewline
150 & -1.1 & 8.6 & 8.8 (8.1) & 8.7 & 94.8\tabularnewline
250 & -1.1 & 6.6 & 6.7 (6.2) & 6.7 & 94.9\tabularnewline
375 & -1.0 & 5.4 & 5.5 (5.1) & 5.4 & 95.0\tabularnewline
500 & -1.1 & 4.6 & 4.7 (4.4) & 4.7 & 95.0\tabularnewline
\end{tabular}

~

~

\begin{tabular}{|cccccc|}
$N_{B(R)}=6,750$ &  &  &  &  & \multicolumn{1}{c}{}\tabularnewline
\cline{1-1} 
Approximate relative bias -1.4 &  &  &  &  & \multicolumn{1}{c}{}\tabularnewline
\cline{1-1} 
$n$ & $RB$ & $ARRSME$ & $RRMSE$ & $ERRMSEE$ & $COV95$\tabularnewline
\hline 
75 & -1.5 & 9.3 & 9.6 (11.5) & 9.6 & 94.7\tabularnewline
100 & -1.4 & 8.1 & 8.3 (9.9) & 8.3 & 94.9\tabularnewline
150 & -1.3 & 6.6 & 6.7 (8.1) & 6.7 & 94.6\tabularnewline
250 & -1.4 & 5.1 & 5.3 (6.2) & 5.1 & 94.3\tabularnewline
375 & -1.3 & 4.1 & 4.3 (5.1) & 4.2 & 94.1\tabularnewline
500 & -1.4 & 3.5 & 3.8 (4.4) & 3.6 & 93.6\tabularnewline
\end{tabular}
\end{table}

~

\begin{table}[H]
\caption{Percentage values of relative bias ($RB$), first order approximation
of relative root mean squared error ($ARRMSE$), relative root mean
squared error ($RRMSE$), expectation of the relative root mean squared
error estimator ($ERRMSEE$) and coverage of the nominal 95\% confidence
intervals achieved from a population of 10,000 units, a sampled sub-population
of 7,500 units with 2,250; 4,500 and 7,500 respondent units, sample
sizes $n=75;100;150;250;375;500$ selected by means of simple random
sampling without replacement and correlations $\rho_{XY}=0.3$, $\rho_{ZY}=0.9$
and $\rho_{XZ}=0$. Values in parentheses are the $RRMSE$s of the
Horvitz-Thompson estimator in absence of nonresponse and under-coverage.}

~

\begin{tabular}{|cccccc|}
$N_{B(R)}=2,250$ &  &  &  &  & \multicolumn{1}{c}{}\tabularnewline
\cline{1-1} 
Approximate relative bias 0.6 &  &  &  &  & \multicolumn{1}{c}{}\tabularnewline
\cline{1-1} 
$n$ & $RB$ & $ARRSME$ & $RRMSE$ & $ERRMSEE$ & $COV95$\tabularnewline
\hline 
75 & 0.6 & 13.9 & 14.7 (11.5) & 14.2 & 94.1\tabularnewline
100 & 0.5 & 12.0 & 12.4 (9.9) & 12.2 & 94.5\tabularnewline
150 & 0.5 & 9.8 & 10.1 (8.1) & 9.8 & 94.2\tabularnewline
250 & 0.6 & 7.5 & 7.6 (6.2) & 7.5 & 95.2\tabularnewline
375 & 0.7 & 6.1 & 6.2 (5.1) & 6.1 & 94.9\tabularnewline
500 & 0.5 & 5.2 & 5.3 (4.4) & 5.2 & 95.2\tabularnewline
\end{tabular}

~

~

\begin{tabular}{|cccccc|}
$N_{B(R)}=4,500$ &  &  &  &  & \multicolumn{1}{c}{}\tabularnewline
\cline{1-1} 
Approximate relative bias -0.9 &  &  &  &  & \multicolumn{1}{c}{}\tabularnewline
\cline{1-1} 
$n$ & $RB$ & $ARRSME$ & $RRMSE$ & $ERRMSEE$ & $COV95$\tabularnewline
\hline 
75 & -0.9 & 8.1 & 8.5 (11.5) & 8.4 & 95.3\tabularnewline
100 & -0.9 & 7.0 & 7.2 (9.9) & 7.2 & 95.2\tabularnewline
150 & -1.0 & 5.7 & 5.9 (8.1) & 5.8 & 95.1\tabularnewline
250 & -0.9 & 4.4 & 4.5 (6.2) & 4.5 & 95.3\tabularnewline
375 & -0.9 & 3.5 & 3.7 (5.1) & 3.6 & 94.6\tabularnewline
500 & -0.9 & 3.0 & 3.2 (4.4) & 3.1 & 94.3\tabularnewline
\end{tabular}

~

~

\begin{tabular}{|cccccc|}
$N_{B(R)}=6,750$ &  &  &  &  & \multicolumn{1}{c}{}\tabularnewline
\cline{1-1} 
Approximate relative bias -0.9 &  &  &  &  & \multicolumn{1}{c}{}\tabularnewline
\cline{1-1} 
$n$ & $RB$ & $ARRSME$ & $RRMSE$ & $ERRMSEE$ & $COV95$\tabularnewline
\hline 
75 & -1.0 & 4.8 & 5.0 (11.5) & 5.0 & 95.2\tabularnewline
100 & -0.9 & 4.1 & 4.3 (9.9) & 4.3 & 95.1\tabularnewline
150 & -0.9 & 3.4 & 3.5 (8.1) & 3.4 & 94.7\tabularnewline
250 & -1.0 & 2.6 & 2.8 (6.2) & 2.6 & 94.1\tabularnewline
375 & -0.9 & 2.1 & 2.3 (5.1) & 2.1 & 93.5\tabularnewline
500 & -0.9 & 1.8 & 2.0 (4.4) & 1.8 & 92.6\tabularnewline
\end{tabular}
\end{table}

~

\begin{table}[H]
\caption{Percentage values of relative bias ($RB$), first order approximation
of relative root mean squared error ($ARRMSE$), relative root mean
squared error ($RRMSE$), expectation of the relative root mean squared
error estimator ($ERRMSEE$) and coverage of the nominal 95\% confidence
intervals achieved from a population of 10,000 units, a sampled sub-population
of 7,500 units with 2,250; 4,500 and 7,500 respondent units, sample
sizes $n=75;100;150;250;375;500$ selected by means of simple random
sampling without replacement and correlations $\rho_{XY}=0.6$, $\rho_{ZY}=0.3$
and $\rho_{XZ}=0$. Values in parentheses are the $RRMSE$s of the
Horvitz-Thompson estimator in absence of nonresponse and under-coverage.}

~

\begin{tabular}{|cccccc|}
$N_{B(R)}=2,250$ &  &  &  &  & \multicolumn{1}{c}{}\tabularnewline
\cline{1-1} 
Approximate relative bias -1.8 &  &  &  &  & \multicolumn{1}{c}{}\tabularnewline
\cline{1-1} 
$n$ & $RB$ & $ARRSME$ & $RRMSE$ & $ERRMSEE$ & $COV95$\tabularnewline
\hline 
75 & -1.8 & 16.8 & 17.7 (11.5) & 17.3 & 92.5\tabularnewline
100 & -1.9 & 14.5 & 14.9 (9.9) & 14.9 & 93.5\tabularnewline
150 & -1.7 & 11.8 & 12.2 (8.1) & 12.1 & 93.6\tabularnewline
250 & -1.8 & 9.1 & 9.4 (6.2) & 9.3 & 94.6\tabularnewline
375 & -1.7 & 7.4 & 7.6 (5.1) & 7.5 & 94.4\tabularnewline
500 & -1.9 & 6.3 & 6.6 (4.4) & 6.4 & 94.1\tabularnewline
\end{tabular}

~

~

\begin{tabular}{|cccccc|}
$N_{B(R)}=4,500$ &  &  &  &  & \multicolumn{1}{c}{}\tabularnewline
\cline{1-1} 
Approximate relative bias -1.3 &  &  &  &  & \multicolumn{1}{c}{}\tabularnewline
\cline{1-1} 
$n$ & $RB$ & $ARRSME$ & $RRMSE$ & $ERRMSEE$ & $COV95$\tabularnewline
\hline 
75 & -1.3 & 11.8 & 12.1 (11.5) & 12.1 & 94.2\tabularnewline
100 & -1.2 & 10.2 & 10.2 (9.9) & 10.4 & 94.8\tabularnewline
150 & -1.3 & 8.3 & 8.5 (8.1) & 8.5 & 94.6\tabularnewline
250 & -1.3 & 6.4 & 6.5 (6.2) & 6.5 & 94.6\tabularnewline
375 & -1.3 & 5.2 & 5.4 (5.1) & 5.2 & 94.6\tabularnewline
500 & -1.3 & 4.4 & 4.6 (4.4) & 4.5 & 94.7\tabularnewline
\end{tabular}

~

~

\begin{tabular}{|cccccc|}
$N_{B(R)}=6,750$ &  &  &  &  & \multicolumn{1}{c}{}\tabularnewline
\cline{1-1} 
Approximate relative bias -1.6 &  &  &  &  & \multicolumn{1}{c}{}\tabularnewline
\cline{1-1} 
$n$ & $RB$ & $ARRSME$ & $RRMSE$ & $ERRMSEE$ & $COV95$\tabularnewline
\hline 
75 & -1.7 & 9.4 & 9.6 (11.5) & 9.6 & 94.4\tabularnewline
100 & -1.6 & 8.1 & 8.3 (9.9) & 8.3 & 94.6\tabularnewline
150 & -1.6 & 6.6 & 6.8 (8.1) & 6.7 & 94.7\tabularnewline
250 & -1.7 & 5.1 & 5.4 (6.2) & 5.2 & 93.9\tabularnewline
375 & -1.6 & 4.1 & 4.4 (5.1) & 4.2 & 93.7\tabularnewline
500 & -1.6 & 3.5 & 3.9 (4.4) & 3.6 & 92.8\tabularnewline
\end{tabular}
\end{table}

~

\begin{table}[H]
\caption{Percentage values of relative bias ($RB$), first order approximation
of relative root mean squared error ($ARRMSE$), relative root mean
squared error ($RRMSE$), expectation of the relative root mean squared
error estimator ($ERRMSEE$) and coverage of the nominal 95\% confidence
intervals achieved from a population of 10,000 units, a sampled sub-population
of 7,500 units with 2,250; 4,500 and 7,500 respondent units, sample
sizes $n=75;100;150;250;375;500$ selected by means of simple random
sampling without replacement and correlations $\rho_{XY}=0.6$, $\rho_{ZY}=0.6$
and $\rho_{XZ}=0$. Values in parentheses are the $RRMSE$s of the
Horvitz-Thompson estimator in absence of nonresponse and under-coverage.}

~

\begin{tabular}{|cccccc|}
$N_{B(R)}=2,250$ &  &  &  &  & \multicolumn{1}{c}{}\tabularnewline
\cline{1-1} 
Approximate relative bias -0.8 &  &  &  &  & \multicolumn{1}{c}{}\tabularnewline
\cline{1-1} 
$n$ & $RB$ & $ARRSME$ & $RRMSE$ & $ERRMSEE$ & $COV95$\tabularnewline
\hline 
75 & -1.7 & 9.4 & 9.6 (11.5) & 9.6 & 94.4\tabularnewline
100 & -1.6 & 8.1 & 8.3 (9.9) & 8.3 & 94.6\tabularnewline
150 & -1.6 & 6.6 & 6.8 (8.1) & 6.7 & 94.7\tabularnewline
250 & -1.7 & 5.1 & 5.4 (6.2) & 5.2 & 93.9\tabularnewline
375 & -1.6 & 4.1 & 4.4 (5.1) & 4.2 & 93.7\tabularnewline
500 & -1.6 & 3.5 & 3.9 (4.4) & 3.6 & 92.8\tabularnewline
\end{tabular}

~

~

\begin{tabular}{|cccccc|}
$N_{B(R)}=4,500$ &  &  &  &  & \multicolumn{1}{c}{}\tabularnewline
\cline{1-1} 
Approximate relative bias -1.3 &  &  &  &  & \multicolumn{1}{c}{}\tabularnewline
\cline{1-1} 
$n$ & $RB$ & $ARRSME$ & $RRMSE$ & $ERRMSEE$ & $COV95$\tabularnewline
\hline 
75 & -1.3 & 9.5 & 9.8 (11.5) & 9.7 & 94.3\tabularnewline
100 & -1.2 & 8.2 & 8.3 (9.9) & 8.4 & 95.1\tabularnewline
150 & -1.3 & 6.7 & 6.9 (8.1) & 6.8 & 94.5\tabularnewline
250 & -1.3 & 5.1 & 5.3 (6.2) & 5.2 & 94.7\tabularnewline
375 & -1.3 & 4.1 & 4.4 (5.1) & 4.2 & 94.3\tabularnewline
500 & -1.3 & 3.6 & 3.8 (4.4) & 3.6 & 93.9\tabularnewline
\end{tabular}

~

~

\begin{tabular}{|cccccc|}
$N_{B(R)}=6,750$ &  &  &  &  & \multicolumn{1}{c}{}\tabularnewline
\cline{1-1} 
Approximate relative bias -1.4 &  &  &  &  & \multicolumn{1}{c}{}\tabularnewline
\cline{1-1} 
$n$ & $RB$ & $ARRSME$ & $RRMSE$ & $ERRMSEE$ & $COV95$\tabularnewline
\hline 
75 & -1.5 & 7.0 & 7.3 (11.5) & 7.2 & 94.7\tabularnewline
100 & -1.4 & 6.0 & 6.3 (9.9) & 6.2 & 94.7\tabularnewline
150 & -1.4 & 4.9 & 5.2 (8.1) & 5.0 & 94.7\tabularnewline
250 & -1.5 & 3.8 & 4.1 (6.2) & 3.9 & 93.8\tabularnewline
375 & -1.4 & 3.1 & 3.4 (5.1) & 3.1 & 93.3\tabularnewline
500 & -1.4 & 2.6 & 3.0 (4.4) & 2.7 & 91.8\tabularnewline
\end{tabular}
\end{table}

~

\begin{table}[H]
\caption{Percentage values of relative bias ($RB$), first order approximation
of relative root mean squared error ($ARRMSE$), relative root mean
squared error ($RRMSE$), expectation of the relative root mean squared
error estimator ($ERRMSEE$) and coverage of the nominal 95\% confidence
intervals achieved from a population of 10,000 units, a sampled sub-population
of 7,500 units with 2,250; 4,500 and 7,500 respondent units, sample
sizes $n=75;100;150;250;375;500$ selected by means of simple random
sampling without replacement and correlations $\rho_{XY}=0.6$, $\rho_{ZY}=0.9$
and $\rho_{XZ}=0$. Values in parentheses are the $RRMSE$s of the
Horvitz-Thompson estimator in absence of nonresponse and under-coverage.}

~

\begin{tabular}{|cccccc|}
$N_{B(R)}=2,250$ &  &  &  &  & \multicolumn{1}{c}{}\tabularnewline
\cline{1-1} 
Approximate relative bias 0.6 &  &  &  &  & \multicolumn{1}{c}{}\tabularnewline
\cline{1-1} 
$n$ & $RB$ & $ARRSME$ & $RRMSE$ & $ERRMSEE$ & $COV95$\tabularnewline
\hline 
75 & 0.6 & 10.5 & 11.4 (11.5) & 11.0 & 94.6\tabularnewline
100 & 0.6 & 9.1 & 9.6 (9.9) & 9.4 & 94.9\tabularnewline
150 & 0.5 & 7.4 & 7.7 (8.1) & 7.6 & 95.1\tabularnewline
250 & 0.6 & 5.7 & 5.8 (6.2) & 5.8 & 95.3\tabularnewline
375 & 0.6 & 4.6 & 4.8 (5.1) & 4.6 & 95.2\tabularnewline
500 & 0.5 & 4.0 & 4.1 (4.4) & 4.0 & 95.0\tabularnewline
\end{tabular}

~

~

\begin{tabular}{|cccccc|}
$N_{B(R)}=4,500$ &  &  &  &  & \multicolumn{1}{c}{}\tabularnewline
\cline{1-1} 
Approximate relative bias -1.0 &  &  &  &  & \multicolumn{1}{c}{}\tabularnewline
\cline{1-1} 
$n$ & $RB$ & $ARRSME$ & $RRMSE$ & $ERRMSEE$ & $COV95$\tabularnewline
\hline 
75 & -1.0 & 5.7 & 6.1 (11.5) & 6.2 & 95.8\tabularnewline
100 & -1.0 & 5.0 & 5.3 (9.9) & 5.3 & 95.3\tabularnewline
150 & -1.0 & 4.0 & 4.3 (8.1) & 4.2 & 95.2\tabularnewline
250 & -1.0 & 3.1 & 3.3 (6.2) & 3.2 & 94.8\tabularnewline
375 & -1.0 & 2.5 & 2.7 (5.1) & 2.6 & 93.8\tabularnewline
500 & -1.0 & 2.2 & 2.4 (4.4) & 2.2 & 93.2\tabularnewline
\end{tabular}

~

~

\begin{tabular}{|cccccc|}
$N_{B(R)}=6,750$ &  &  &  &  & \multicolumn{1}{c}{}\tabularnewline
\cline{1-1} 
Approximate relative bias -0.8 &  &  &  &  & \multicolumn{1}{c}{}\tabularnewline
\cline{1-1} 
$n$ & $RB$ & $ARRSME$ & $RRMSE$ & $ERRMSEE$ & $COV95$\tabularnewline
\hline 
75 & -0.8 & 2.6 & 3.0 (11.5) & 3.0 & 96.6\tabularnewline
100 & -0.8 & 2.3 & 2.5 (9.9) & 2.5 & 96.1\tabularnewline
150 & -0.8 & 1.8 & 2.1 (8.1) & 2.0 & 94.7\tabularnewline
250 & -0.8 & 1.4 & 1.7 (6.2) & 1.5 & 92.7\tabularnewline
375 & -0.8 & 1.1 & 1.4 (5.1) & 1.2 & 90.8\tabularnewline
500 & -0.8 & 1.0 & 1.3 (4.4) & 1.0 & 88.5\tabularnewline
\end{tabular}
\end{table}

~

\begin{table}[H]
\caption{Percentage values of relative bias ($RB$), first order approximation
of relative root mean squared error ($ARRMSE$), relative root mean
squared error ($RRMSE$), expectation of the relative root mean squared
error estimator ($ERRMSEE$) and coverage of the nominal 95\% confidence
intervals achieved from a population of 10,000 units, a sampled sub-population
of 7,500 units with 2,250; 4,500 and 7,500 respondent units, sample
sizes $n=75;100;150;250;375;500$ selected by means of simple random
sampling without replacement and correlations $\rho_{XY}=0.9$, $\rho_{ZY}=0.3$
and $\rho_{XZ}=0$. Values in parentheses are the $RRMSE$s of the
Horvitz-Thompson estimator in absence of nonresponse and under-coverage.}

~

\begin{tabular}{|cccccc|}
$N_{B(R)}=2,250$ &  &  &  &  & \multicolumn{1}{c}{}\tabularnewline
\cline{1-1} 
Approximate relative bias -1.3 &  &  &  &  & \multicolumn{1}{c}{}\tabularnewline
\cline{1-1} 
$n$ & $RB$ & $ARRSME$ & $RRMSE$ & $ERRMSEE$ & $COV95$\tabularnewline
\hline 
75 & -1.3 & 9.3 & 10.1 (11.5) & 9.7 & 93.9\tabularnewline
100 & -1.3 & 8.0 & 8.5 (9.9) & 8.3 & 94.4\tabularnewline
150 & -1.3 & 6.5 & 6.9 (8.1) & 6.7 & 94.1\tabularnewline
250 & -1.3 & 5.0 & 5.3 (6.2) & 5.1 & 94.2\tabularnewline
375 & -1.3 & 4.1 & 4.3 (5.1) & 4.1 & 94.4\tabularnewline
500 & -1.4 & 3.5 & 3.8 (4.4) & 3.6 & 93.3\tabularnewline
\end{tabular}

~

~

\begin{tabular}{|cccccc|}
$N_{B(R)}=4,500$ &  &  &  &  & \multicolumn{1}{c}{}\tabularnewline
\cline{1-1} 
Approximate relative bias -1.4 &  &  &  &  & \multicolumn{1}{c}{}\tabularnewline
\cline{1-1} 
$n$ & $RB$ & $ARRSME$ & $RRMSE$ & $ERRMSEE$ & $COV95$\tabularnewline
\hline 
75 & -1.4 & 6.3 & 6.8 (11.5) & 6.7 & 94.4\tabularnewline
100 & -1.4 & 5.5 & 5.7 (9.9) & 5.7 & 94.6\tabularnewline
150 & -1.4 & 4.4 & 4.7 (8.1) & 4.6 & 94.1\tabularnewline
250 & -1.4 & 3.4 & 3.7 (6.2) & 3.5 & 93.7\tabularnewline
375 & -1.4 & 2.8 & 3.1 (5.1) & 2.8 & 92.7\tabularnewline
500 & -1.4 & 2.4 & 2.8 (4.4) & 2.4 & 91.4\tabularnewline
\end{tabular}

~

~

\begin{tabular}{|cccccc|}
$N_{B(R)}=6,750$ &  &  &  &  & \multicolumn{1}{c}{}\tabularnewline
\cline{1-1} 
Approximate relative bias -1.4 &  &  &  &  & \multicolumn{1}{c}{}\tabularnewline
\cline{1-1} 
$n$ & $RB$ & $ARRSME$ & $RRMSE$ & $ERRMSEE$ & $COV95$\tabularnewline
\hline 
75 & -1.4 & 4.9 & 5.2 (11.5) & 5.2 & 94.5\tabularnewline
100 & -1.4 & 4.2 & 4.5 (9.9) & 4.4 & 94.4\tabularnewline
150 & -1.4 & 3.4 & 3.7 (8.1) & 3.5 & 93.8\tabularnewline
250 & -1.4 & 2.6 & 3.0 (6.2) & 2.7 & 92.5\tabularnewline
375 & -1.3 & 2.1 & 2.5 (5.1) & 2.2 & 91.2\tabularnewline
500 & -1.4 & 1.8 & 2.3 (4.4) & 1.9 & 88.4\tabularnewline
\end{tabular}
\end{table}

~

\begin{table}[H]
\caption{Percentage values of relative bias ($RB$), first order approximation
of relative root mean squared error ($ARRMSE$), relative root mean
squared error ($RRMSE$), expectation of the relative root mean squared
error estimator ($ERRMSEE$) and coverage of the nominal 95\% confidence
intervals achieved from a population of 10,000 units, a sampled sub-population
of 7,500 units with 2,250; 4,500 and 7,500 respondent units, sample
sizes $n=75;100;150;250;375;500$ selected by means of simple random
sampling without replacement and correlations $\rho_{XY}=0.9$, $\rho_{ZY}=0.6$
and $\rho_{XZ}=0$. Values in parentheses are the $RRMSE$s of the
Horvitz-Thompson estimator in absence of nonresponse and under-coverage.}

~

\begin{tabular}{|cccccc|}
$N_{B(R)}=2,250$ &  &  &  &  & \multicolumn{1}{c}{}\tabularnewline
\cline{1-1} 
Approximate relative bias -0.7 &  &  &  &  & \multicolumn{1}{c}{}\tabularnewline
\cline{1-1} 
$n$ & $RB$ & $ARRSME$ & $RRMSE$ & $ERRMSEE$ & $COV95$\tabularnewline
\hline 
75 & -0.7 & 6.7 & 7.6 (11.5) & 7.5 & 95.4\tabularnewline
100 & -0.7 & 5.8 & 6.3 (9.9) & 6.3 & 95.6\tabularnewline
6150 & -0.8 & 4.7 & 5.0 (8.1) & 5.0 & 95.3\tabularnewline
250 & -0.7 & 3.6 & 3.8 (6.2) & 3.8 & 95.1\tabularnewline
375 & -0.7 & 2.9 & 3.1 (5.1) & 3.0 & 95.2\tabularnewline
500 & -0.7 & 2.5 & 2.7 (4.4) & 2.6 & 94.1\tabularnewline
\end{tabular}

~

~

\begin{tabular}{|cccccc|}
$N_{B(R)}=4,500$ &  &  &  &  & \multicolumn{1}{c}{}\tabularnewline
\cline{1-1} 
Approximate relative bias -1.3 &  &  &  &  & \multicolumn{1}{c}{}\tabularnewline
\cline{1-1} 
$n$ & $RB$ & $ARRSME$ & $RRMSE$ & $ERRMSEE$ & $COV95$\tabularnewline
\hline 
75 & -1.4 & 4.1 & 4.7 (11.5) & 4.7 & 94.9\tabularnewline
100 & -1.3 & 3.5 & 4.0 (9.9) & 3.9 & 95.0\tabularnewline
150 & -1.4 & 2.9 & 3.3 (8.1) & 3.1 & 94.0\tabularnewline
250 & -1.3 & 2.2 & 2.6 (6.2) & 2.3 & 92.2\tabularnewline
375 & -1.4 & 1.8 & 2.3 (5.1) & 1.9 & 89.5\tabularnewline
500 & -1.3 & 1.5 & 2.1 (4.4) & 1.6 & 86.8\tabularnewline
\end{tabular}

~

~

\begin{tabular}{|cccccc|}
$N_{B(R)}=6,750$ &  &  &  &  & \multicolumn{1}{c}{}\tabularnewline
\cline{1-1} 
Approximate relative bias -1.9 &  &  &  &  & \multicolumn{1}{c}{}\tabularnewline
\cline{1-1} 
$n$ & $RB$ & $ARRSME$ & $RRMSE$ & $ERRMSEE$ & $COV95$\tabularnewline
\hline 
75 & -1.1 & 2.7 & 3.2 (11.5) & 3.2 & 95.9\tabularnewline
100 & -1.2 & 2.4 & 2.8 (9.9) & 2.7 & 94.6\tabularnewline
150 & -1.2 & 1.9 & 2.3 (8.1) & 2.1 & 92.5\tabularnewline
250 & -1.2 & 1.5 & 1.9 (6.2) & 1.6 & 89.2\tabularnewline
375 & -1.2 & 1.2 & 1.7 (5.1) & 1.2 & 85.9\tabularnewline
500 & -1.2 & 1.0 & 1.6 (4.4) & 1.1 & 81.2\tabularnewline
\end{tabular}
\end{table}

~

\begin{table}[H]
\caption{Percentage values of relative bias ($RB$), first order approximation
of relative root mean squared error ($ARRMSE$), relative root mean
squared error ($RRMSE$), expectation of the relative root mean squared
error estimator ($ERRMSEE$) and coverage of the nominal 95\% confidence
intervals achieved from a population of 10,000 units, a sampled sub-population
of 7,500 units with 2,250; 4,500 and 7,500 respondent units, sample
sizes $n=75;100;150;250;375;500$ selected by means of simple random
sampling without replacement and correlations $\rho_{XY}=0.9$, $\rho_{ZY}=0.9$
and $\rho_{XZ}=0$. Values in parentheses are the $RRMSE$s of the
Horvitz-Thompson estimator in absence of nonresponse and under-coverage.}

~

\begin{tabular}{|cccccc|}
$N_{B(R)}=2,250$ &  &  &  &  & \multicolumn{1}{c}{}\tabularnewline
\cline{1-1} 
Approximate relative bias -0.6 &  &  &  &  & \multicolumn{1}{c}{}\tabularnewline
\cline{1-1} 
$n$ & $RB$ & $ARRSME$ & $RRMSE$ & $ERRMSEE$ & $COV95$\tabularnewline
\hline 
75 & -0.6 & 8.3 & 8.7 (11.5) & 8.7 & 95.3\tabularnewline
100 & -0.7 & 7.2 & 7.4 (9.9) & 7.5 & 95.5\tabularnewline
150 & -0.6 & 5.8 & 6.0 (8.1) & 6.0 & 95.0\tabularnewline
250 & -0.6 & 4.5 & 4.6 (6.2) & 4.6 & 95.3\tabularnewline
375 & -0.6 & 3.6 & 3.7 (5.1) & 3.7 & 94.9\tabularnewline
500 & -0.7 & 3.1 & 3.2 (4.4) & 3.2 & 95.1\tabularnewline
\end{tabular}

~

~

\begin{tabular}{|cccccc|}
$N_{B(R)}=4,500$ &  &  &  &  & \multicolumn{1}{c}{}\tabularnewline
\cline{1-1} 
Approximate relative bias -1.0 &  &  &  &  & \multicolumn{1}{c}{}\tabularnewline
\cline{1-1} 
$n$ & $RB$ & $ARRSME$ & $RRMSE$ & $ERRMSEE$ & $COV95$\tabularnewline
\hline 
75 & -0.9 & 6.8 & 7.0 (11.5) & 7.0 & 95.4\tabularnewline
100 & -0.9 & 5.9 & 6.0 (9.9) & 6.0 & 95.2\tabularnewline
150 & -1.0 & 4.8 & 4.9 (8.1) & 4.9 & 95.2\tabularnewline
250 & -0.9 & 3.7 & 3.8 (6.2) & 3.7 & 95.1\tabularnewline
375 & -1.0 & 3.0 & 3.1 (5.1) & 3.0 & 94.2\tabularnewline
500 & -0.9 & 2.6 & 2.7 (4.4) & 2.6 & 94.1\tabularnewline
\end{tabular}

~

~

\begin{tabular}{|cccccc|}
$N_{B(R)}=6,750$ &  &  &  &  & \multicolumn{1}{c}{}\tabularnewline
\cline{1-1} 
Approximate relative bias -1.0 &  &  &  &  & \multicolumn{1}{c}{}\tabularnewline
\cline{1-1} 
$n$ & $RB$ & $ARRSME$ & $RRMSE$ & $ERRMSEE$ & $COV95$\tabularnewline
\hline 
75 & -1.0 & 6.2 & 6.4 (11.5) & 6.4 & 95.7\tabularnewline
100 & -1.0 & 5.4 & 5.5 (9.9) & 5.5 & 95.1\tabularnewline
150 & -0.9 & 4.4 & 4.5 (8.1) & 4.4 & 95.1\tabularnewline
250 & -0.9 & 3.4 & 3.5 (6.2) & 3.4 & 94.3\tabularnewline
375 & -1.0 & 2.7 & 2.9 (5.1) & 2.8 & 94.3\tabularnewline
500 & -0.9 & 2.3 & 2.5 (4.4) & 2.4 & 93.6\tabularnewline
\end{tabular}
\end{table}

\textcolor{black}{\pagebreak}

\textcolor{black}{}
\textcolor{black}{\bibliographystyle{apalike}
\bibliography{prova1}
}